\documentclass[a5paper, 10pt]{article}

\usepackage{cite, amsmath, amssymb, booktabs, lastpage, graphicx}
\usepackage[hidelinks]{hyperref}

\usepackage{geometry}
\usepackage{setspace}

\usepackage{amsthm}

\theoremstyle{plain}
\newtheorem{lemma}{Lemma}
\newtheorem{theorem}{Theorem}

\newtheorem{proposition}[theorem]{Proposition}
\newtheorem*{corollary}{Corollary}

\theoremstyle{remark}

\theoremstyle{definition}

\usepackage{graphicx}
\makeatletter

\newcommand{\address}[1]{\def\@address{#1}}
\newcommand{\email}[1]{\def\@email{#1}}

\makeatother

\newcommand{\acknowledgment}[1]{\vspace{5mm}\singlespacing
	{\noindent\textbf{\textit{Acknowledgment\/}:} #1}
}

\usepackage{fancyhdr}

\usepackage[margin=1cm,%
			font=footnotesize,%
			format=hang,%
			labelsep=period,%
			labelfont=bf]{caption}

\title{Minimizing the Sombor Index among Trees with Fixed Degree Sequence}
\author{Mirza Redžić}

\address{$^a$University of Primorska\\

    }

\email{mirza233@gmail.com}

\date{November 13, 2022}

\begin{document}

\maketitle

\begin{abstract}
Vertex-degree-based topological indices have recently gained a lot of attention from mathematical chemists.
One such index that we focus on in this paper is called Sombor index.
After its definition in late 2020, the Sombor index was quickly recognized as a valuable research topic.
In this paper we partially answer the open question of finding the extremal trees with respect to this index for a fixed degree sequence.
Particularly we focus on the lower bound and proceed to show that greedy tree minimizes the Sombor index for a given degree sequence.
\end{abstract}

\onehalfspacing

\section{Introduction}
In 2021, motivated by a geometric approach for interpreting degree-based
topological indices, Gutman introduced a so called \emph{Sombor index} \cite{gutman2021geometric}.

Given a graph $G = (V,E)$, Sombor index of $G$ was defined as 
$${\rm SO}(G) = \sum_{v_iv_j\in E}\sqrt{d(v_i)^2 + d(v_j)^2}$$
This topological index quickly gained a lot of attention from mathematical chemists and most of its mathematical properties have been studied (see e.g. \cite{das2021sombor,cruz2021sombor,ghanbari2021sombor,gutman2021some,reti2021sombor,rada2021general,milovanovic2021some,liu2022sombor}). 
Also, different applications of Sombor index in chemistry have been established \cite{liu2021more,redvzepovic2021chemical}.

In this paper we focus on finding a graph that minimizes the Sombor index among trees with given degree sequence $(d_1,\dots, d_k)$. We show that given a degree sequence $\mathcal{D} = (d_1,\dots, d_k)$, greedy tree minimizes Sombor index among trees with degree sequence $\mathcal{D}$.

\section{Preliminaries and Definitions}
Before proceeding to the main results, we first overview some basic definitions and notations used in the paper.

Note that unless stated otherwise, we are only considering the finite, simple, connected, undirected graphs and given a graph $G$, we denote the vertex set and the edge set of $G$ by $V(G)$ and $E(G)$ respectively.
Two vertices $u$ and $v$ are \emph{adjacent} in a graph $G$ if there is an edge between them. 
The \emph{neighbourhood} of a vertex $v$, denoted $N(v)$, is the set of all vertices that are adjacent to $v$. The \emph{degree} of a vertex $v$, denoted $\deg(v)$, is equal to the cardinality of its neighbourhood. The \emph{closed neighbourhood} of a vertex $v$, denoted $N[v]$, is the neighbourhood of $v$ together with $v$ itself.
For a vertex $v$, if $d(v) = 1$, we say that $v$ is \emph{pendant} vertex.
For a pair of vertices $u,v\in V(G)$, \emph{distance} between $u$ and $v$ is the length of a shortest path between $u$ and $v$ in $G$
Let $L_i$ denote the set of vertices whose minimal distance to a pendant vertex is $i$. 

If $G$ is a graph, we say that $G$ has the \emph{degree sequence} $(d_1,\dots,d_k)$ if for any $1\leq i \leq n$ we can find a vertex in $G$ $v_i$ such that $d(v_i) = d_i$ and $G$ has exactly $k$ non-pendant vertices. 
By convention, we always write the degree sequence in a decreasing order and omit the pendant vertices.
 
Let $G_1,\dots, G_n$ be a sequence of graphs. 
If for some $1\leq i\leq n$, $v$ is a pendant vertex in $G_i$ and non-pendant in $G_{i+1}$, we say that $v$ is \emph{promoted} in $G_{i+1}$.

\section{Main Results}
Let $f(x,y) = \sqrt{x^2 + y^2}$ be defined on positive integers and define 
$$g_{a,b}(x) = f(x,a)-f(x,b)$$
for fixed $b>a\geq 1$. 

\begin{lemma}\label{lemma:g-increasing}
For any fixed positive integers $a,b$, such that $b>a$, it holds that $g_{a,b}$ is a strictly increasing function on its domain.
\end{lemma}
\begin{proof}
We compute 
$$\frac{\partial}{\partial x}g_{a,b}(x) = x\big (\frac{1}{\sqrt{a^2+x^2}} - \frac{1}{\sqrt{b^2+x^2}} \big )$$
and since $x$ is positive and ${\sqrt{a^2+x^2}}<{\sqrt{b^2+x^2}} \Rightarrow \frac{1}{\sqrt{a^2+x^2}} - \frac{1}{\sqrt{b^2+x^2}}>0$, it holds that $\frac{\partial}{\partial x}g_{a,b}(x)>0$. To conclude the proof, it suffices to observe that $g_{a,b}$ is injective on its domain. 
\end{proof}

\begin{proposition}\label{prop:path-degrees}
Let $\mathcal D$ be a fixed sequence of degrees and $T$ be a tree that minimizes the Sombor index among the trees with degree sequence $\mathcal{D}$. 
Let $P = v_1v_2\dots v_t$ be a path in $T$ for some $t\geq 4$ and assume $d(v_1)<d(v_t)$. Then $d(v_2)\leq d(v_{t-1})$.
\end{proposition}
\begin{proof}
Assume for contradiction that $d(v_2)>d(v_{t-1})$.

Let $T'$ be a tree obtained from $T$ by removing the edges $v_1v_2$ and $v_{t-1}v_t$ and adding the edges $v_1v_{t-1}$ and $v_2v_t$. 
It is clear that the degree sequence of $T'$ is $\mathcal D$.

Consider the difference in Sombor indices of $T$ and $T'$.

$${\rm SO}(T) - {\rm SO}(T') = \sum_{uv \in E(T)}f(d(u),d(v)) - \sum_{u'v' \in E(T')}f(d(u'),d(v'))$$
$$= \Big ( f\big (d(v_1),d(v_2)\big )+f\big (d(v_{t-1}),d(v_t)\big ) \Big )-\Big (f\big (d(v_1),d(v_{t-1})\big )+f\big (d(v_2),d(v_t)\big )\Big )$$
$$ = \Big( f\big (d(v_2), d(v_1)\big )- f \big(d(v_2),d(v_t)\big ) \Big ) - \Big( f\big (d(v_{t-1}),d(v_{1})\big ) - f\big (d(v_{t-1}),d(v_{t})\big ) \Big)$$
$$ = g_{d(v_1),d(v_t)}\big(d(v_2)\big)- g_{d(v_1),d(v_t)}\big(d(v_{t-1})\big)>0$$
where the last inequality follows from Lemma \ref{lemma:g-increasing} and the second to last equality follows by rearranging terms and noticing that $f(x,y) = f(y,x)$. 

The proof is concluded by observing that the obtained inequality contradicts the minimality of $T$. 
\end{proof} 

If $T$ is a tree such that for any path $P=v_1\dots v_t$ for some $t\geq 4$ with $d(v_1)<d(v_t)$ it holds that $d(v_2)\leq d(v_{t-1})$, Delorme et. al. proved in \cite{trees} that the following lemmas hold. 

\begin{lemma} 
Let $T$ be as above. There is no path $v_1\dots v_t$ for $t\geq 3$ such that $d(v_1),d(v_t)>d(v_i)$ for some $i\in\{2,\dots, t-1\}$.
\end{lemma}

\begin{lemma}\label{lemma:subtree}
Let $T$ be as above. For every positive integer $d$, the vertices of $T$ with degree at least $d$ form a subtree of $T$.
\end{lemma}

\begin{lemma}
Let $T$ be as above. There are no two non-incident vertices $v_1v_2$ and $v_3v_4$ such that $d(v_1)<d(v_3)\leq d(v_4)<d(v_2)$.
\end{lemma}

Now, by Proposition \ref{prop:path-degrees}, it is clear that all of the previous lemmas hold for a tree that minimizes the Sombor index among the trees with degree sequence $\mathcal D$. 
Furthermore, by Lemma \ref{lemma:subtree}, it is clear that we can root a tree $T$ that minimizes Sombor index among the trees with degree sequence $\mathcal D$, so that for each $i<j$ if $u\in L_i$ and $v\in L_j$, $d(u)\leq d(v)$. 

In \cite{wang}, Wang defines the \emph{greedy tree} as a tree that is achieved by the following algorithm
\begin{itemize}
    \item Label the vertex with the largest degree as $v$ (the root).
    \item Label the neighbours of $v$ as $v_1, v_2, \dots v_{d(v)}$, assign the largest degrees available to them such that $d(v1) \geq d(v2) \geq \dots \geq d(v_{d(v)})$.
    \item Label the neighbours of $v_1$ (except $v$) as $v_{1,1}, v_{1,2},\dots, v_{1,d(v_1)-1}$, so that they take all the largest degrees available and that $d(v11) \geq d(v12) \geq \dots \geq d(v_{1,d(v_1)-1})$, then do the same for $v_2, v_3,\dots v_{d(v)}$.
    \item Repeat the last step for all the newly labelled vertices, always starting with the neighbours of the labelled vertex with largest degree whose neighbours are not yet labelled.
\end{itemize}

We would like to prove that the greedy tree indeed minimizes the Sombor index among trees with degree sequence $\mathcal D$. To do so, we will first prove a few lemmas.

For the rest of this paper, let $d_i$ denote the degree of vertex $v_i$ and $d_i\geq d_j$ whenever $i\leq j$. 
Furthermore, sort the vertices so that if $i<j$ and $d_i = d_j$, $v_i$ denotes the vertex that would be preferred by the greedy algorithm above. 

\begin{lemma}\label{lemma:pendant}
Let $\mathcal{D} = (d_1,\dots, d_k)$ be a degree sequence with $d_i>1$ for every $1\leq i\leq k$. Then we may always assume that if $T$ is a tree that minimizes Sombor index among trees with the degree sequence $\mathcal D$ and is rooted at $v_1$, then every child of $v_k$ is pendant. 
\end{lemma}
\begin{proof}
Let $T$ be a tree that minimizes Sombor index among trees with the degree sequence $\mathcal D$ with inner vertices $v_1,\dots, v_k$ sorted decreasingly by degrees and rooted at $v_1$.

Let $v_j$ be a child of $v_k$ with the largest degree.
It follows from the definition that $d_j\geq d_k$. 
Assume first that $d_j>d_k$. 
Consider the graph $H$ induced by all the vertices of $T$ that have degree at least $d_j$. 
Clearly, the unique path from $v_j$ to $v_1$ contains $v_k$, which implies that $H$ is disconnected.
However, this contradicts Lemma \ref{lemma:subtree}.

Hence $d_j\leq d_k$. 
If $d_j<d_k$, then $d_j = 1$, since $v_k$ has the smallest degree of all non-pendant vertices and we are done.

Otherwise, if $d_j = d_k$, we may relabel $v_j$ and $v_k$ without loss of generality. 

The proof now follows by a simple inductive argument.
\end{proof}

From this last lemma, we may easily obtain the following corollary.

\begin{corollary}\label{cor:min-formula}
Let $T$ be a tree that minimizes Sombor index among trees with degree sequence $\mathcal{D} = (d_1, \dots, d_k)$. Then $T$ can be obtained from a tree $T'$ with degree sequence $\mathcal{D}' = (d_1,\dots, d_{k-1})$ by adding $d_k-1$ pendant vertices as children to a pendant vertex $v$ of $T'$. 
Furthermore, $${\rm SO}(T) = {\rm SO}(T') + (d_k-1)\sqrt{d_k^2 + 1} + \sqrt{d_k^2 + d_p^2}- \sqrt{d_p^2 + 1}$$
where $v_p$ is the parent of $v$ in $T'$.
\end{corollary}
\begin{proof}
The first part follows directly from Lemma \ref{lemma:pendant}.

Consider ${\rm SO}(T) - {\rm SO}(T')$. 
This expression contains $(d_k-1)\sqrt{d_k^2 + 1}$, since we added $d_k-1$ new pendant vertices in $T$. 
Other than that, all the vertices are identical, hence all the terms cancel out, except the pendant child of $v_p$ which was promoted to the inner vertex $v_k$. Hence the edge $v_pv_k$ contributes to $\sqrt{d_p^2 + 1}$ in ${\rm SO}(T')$ and $\sqrt{d_p^2 + d_k^2}$ in ${\rm SO}(T)$. The formula follows.
\end{proof}

\begin{lemma}\label{lemma:h-decreasing}
Let $a>1$ be a positive integer. 
Define a function $h_a$ on positive integers as $h_a(x) = \sqrt{a^2 + x^2}- \sqrt{x^2 + 1}$. 
Then $h_a$ is a strictly decreasing function.
\end{lemma}
\begin{proof}
We notice that $h_a(x) = -g_{1,a}(x)$ and by Lemma \ref{lemma:g-increasing} $g$ is strictly increasing.
\end{proof}

\begin{lemma} \label{lemma:recursion-base}
Let $T$ be a tree with degree sequence $\mathcal{D} = (d_1,\dots, d_k)$ rooted at $v_1$ satisfying the condition of Proposition \ref{prop:path-degrees}. By Lemma \ref{lemma:pendant} we may assume that all the children of $v_k$ are pendant.
Define $T'$ as a tree obtained from $T$ by removing the pendant vertices adjacent to $v_k$. Then $T'$ also satisfies the condition of Proposition \ref{prop:path-degrees}.
\end{lemma}
\begin{proof}
Assume for contradiction that there exists a path $P = (u_1,\dots, u_t)$ in $T'$ with $d(u_1)<d(u_t)$ and $d(u_2)>d(u_{t-1})$. 

If $v_k$ is not in $P$, we are done, since any such path in $T'$ is also contained in $T$. 
Hence, since $v_k$ is pendant in $T'$ by construction, we may assume that $v_k = u_1$.

Since $T$ satisfies the condition of Proposition \ref{prop:path-degrees}, by adding $d_k - 1$ pendant children to $v_k$ in $T'$, it must hold that $d(u_1) = d(u_t)$ in $T$. 
But since we sorted the vertices in the order they would be selected by our greedy algorithm (see the paragraph after the algorithm description), we reach a contradiction, since $v_k$'s parent has the larger degree and would hence be selected first by the algorithm.
\end{proof}

Lemma \ref{lemma:recursion-base} gives us a nice recursive way to think about the trees satisfying Proposition \ref{prop:path-degrees}.
Namely, if $T$ is a tree with the degree sequence $(d_1,\dots, d_k)$ satisfying this Proposition, we define $T_k = T$, and we let $T_{n-1}$ be a tree with degree sequence $d_1,\dots, d_{n-1}$ obtained from $T_n$ by removing the pendant children of $v_n$.

We observe that $T_1$ is always isomorphic to $K_{1,d_1}$ and that given the sequence $T_1,\dots, T_n = T$, with $T$ satisfying the condition of Proposition \ref{prop:path-degrees}, by a simple inductive argument applied to the previous lemma, every tree in this sequence satisfies the condition of this proposition.

\begin{theorem}
Greedy tree minimizes the Sombor index among all the trees with the degree sequence $\mathcal{D}$.
\end{theorem}




\begin{proof}
Let $T$ be the tree obtained by the algorithm for constructing greedy tree after $k$ iterations and let $T'$ be any tree with degree sequence $(d_1,\dots,d_k)$.
We would like to show that ${\rm SO}(T)\leq {\rm SO}(T')$.

If $T'$ has a path $P = u_1\dots u_l$ such that $d(u_1)<d(u_t)$ and $d(u_2)>d(u_{t-1})$, then by Proposition \ref{prop:path-degrees} we can obtain a tree with smaller Sombor index as in the proof of Proposition \ref{prop:path-degrees}. 
Hence, we may assume without loss of generality that $T'$ satisfies the condition of Proposition \ref{prop:path-degrees} and we construct the recursive sequences $T_1,\dots, T_k$ and $T'_1,\dots, T'_k$ for $T$ and $T'$ respectively as above.

Let $t$ be the largest value such that $T_{t-1}\cong T'_{t-1}$. Clearly at least one such $t$ exists, since it always holds that $T_1\cong T'_1$. 

We observe that by Corollary \ref{cor:min-formula}
$${\rm SO}(T_t) - {\rm SO}(T'_t) = (\sqrt{d_p^2+d_t^2} - \sqrt{d_p^2+1}) - (\sqrt{d_q^2+d_t^2} - \sqrt{d_q^2+1})$$
where in $T_t$ the vertex $v_t$ is a child of $v_p$ and in $T'_t$, it is a child of $v_q$.

Since $T_{t-1}$ and $T'_{t-1}$ are isomorphic, and the parent of $v_t$ in $T_t$ ($v_p$) was selected by greedy algorithm, it follows that $d_q\leq d_p$.

More generally, we observe that at iteration $n$, the parent of $v_n$ in $T'_n$ can have degree larger than the parent of $v_n$ in $T_n$ only if at some $i<n$ the parent of $v_i$ in $T'$ has degree smaller than the parent of $v_i$ in $T$ and the number of such $i$'s is strictly larger than the number of indices $j<n$ such that the parent of $v_j$ in $T'_j$ has the smaller degree than the parent of $v_j$ in $T_j$. 

Clearly, at the iteration $t-1$ the number of such $i$'s and $j$'s is the same, so we may only consider what happens from iteration $t$ on.

We fix an iteration $n\geq t$. 
Let $I$ denote the set of all indices $i<n$ where the parent of $v_i$ in $T'$ has degree smaller than the parent of $v_i$ in $T$.
We proceed by strong induction on $|I|$ to show that ${\rm SO}(T_n) \leq {\rm SO}(T'_n)$.

If $|I| = 0$, then the assertion is trivially satisfied.
Now assume that for any tree with $|I|\leq k$, the assertion holds.

Consider a tree with $|I| = k+1$. 
Let $r<n$ denote the largest iteration where the parent of $v_{r}$ in $T'_{r}$ has degree $d_q$ and the parent of $v_{r}$ in $T_{r}$ has degree $d_p$ with $d_q<d_p$. 



At iteration $r-1$, by induction hypothesis, we have that ${\rm SO}(T_{r-1})\leq {\rm SO}(T'_{r-1})$.
Hence,  by Lemma \ref{lemma:h-decreasing}, it follows that ${\rm SO}(T_r) - {\rm SO}(T'_r) \leq h_{d_r}(d_p) - h_{d_r}(d_q) <0$.

We observe now that if $v_p$ remains in $L_1$\footnote{Correct statement here is that $v_p$ has one more pendant child in $T'_n$ than in $T_n$ (which is what we need), but if $n$ is sufficiently large, those statements are equivalent.} in $T'_n$, we can apply the induction hypothesis as follows.
Let $T''_n$ be a tree constructed by taking $T'_n$ and deleting $d_r-1$ vertices from the child of $v_q$ and adding them to $v_p$. By the computation above, if follows that ${\rm SO}(T''_n)\leq {\rm SO}(T'_n)$ and by induction hypothesis ${\rm SO}(T_n)\leq {\rm SO}(T''_n)$.

Hence, we may assume that $v_p$ is not in $L_1$ in $T'_n$. 
Therefore, at some step $r+i$, we have that in $T'_{r+i}$ the child of $v_p$, that was pendant in $T'_r$ and non-pendant in $T_r$, was promoted. 
Say that in $T_{r+i}$ the vertex $v_{p'}$ was promoted.

We consider two cases. Namely $d_{q}\leq d_{p'}$ and $d_{q}> d_{p'}$.

Consider first the case when $d_q\leq d_{p'}$.
\begin{equation}\label{ineq:1}
{\rm SO}(T_{r+i}) - {\rm SO}(T'_{r+i}) \leq h_{d_r}(d_p) - h_{d_r}(d_q) + h_{d_{r+i}}(d_{p'}) - h_{d_{r+i}}(d_{p})
\end{equation} 
\begin{equation*}\label{eq:2}
\begin{split}
= (\sqrt{d_r^2+d_p^2}-\sqrt{d_p^2+1})-(\sqrt{d_r^2+d_q^2}-\sqrt{d_q^2+1}) + \\ (\sqrt{d_{r+i}^2+d_{p'}^2}- \sqrt{d_{p'}^2+1}) - (\sqrt{d_{r+i}^2+d_p^2} - \sqrt{d_p^2+1})
\end{split}\end{equation*}
\begin{equation*}\label{eq:3}
\begin{split}
= (\sqrt{d_r^2+d_p^2} +\sqrt{d_q^2+1} + \sqrt{d_{r+i}^2+d_{p'}^2})\\
-(\sqrt{d_r^2+d_q^2} +  \sqrt{d_{p'}^2+1}+ \sqrt{d_{r+i}^2+d_p^2}) 
\end{split}
\end{equation*}

\begin{equation}\label{ineq:4}
\leq (\sqrt{d_r^2+d_p^2} +\sqrt{d_q^2+1})
-(\sqrt{d_r^2+d_q^2} +  \sqrt{d_{p'}^2+1})
\end{equation}
\begin{equation}\label{ineq:5} 
= h_{d_r}(d_{p'}) - h_{d_r}(d_q)\leq 0
\end{equation}
Where  (\ref{ineq:1}) follows from the induction hypothesis, (\ref{ineq:4}) follows by observing that $d_{p'}\leq d_p$, as otherwise the greedy algorithm would select $v_{p'}$ before $v_p$ and (\ref{ineq:5}) follows from Lemma \ref{lemma:h-decreasing}.

Finally, consider the case when $d_q> d_{p'}$. 
Notice that by construction of our greedy algorithm, $d_q$ has no pendant children in $T_{r+i}$, moreover this implies that the child of $v_q$ that was promoted in $T'_r$ and not promoted in $T_r$ was promoted in $T_{r+j}$, and by maximality of $r$ it follows that the parent of $v_{r+j}$, say $v_s$ in $T'_{r+j}$ has degree at least $d_q$. 

Consider the following algorithm. Note that for simplicity we drop the indices in $T, T', h$ in the algorithm, but we will consider them implicit.
\begin{itemize}
    \item When the child of $v_q$ was promoted in $T$, the child of $v_s$ was promoted in $T'$, and degree of $v_s$ in $T'$ is at least $d_q$. Thus the greedy algorithm already promoted this child of $v_s$ in $T$.
    \item Consider the step at which the child of $v_s$ in $T$ from the previous step was promoted. Suppose that at this time a child of vertex $v_l$ was promoted. If $d_{\ell}\leq d_q$, we stop the algorithm and return $[v_q, v_s, v_l]$ and the indices $[r+j, \_]$, where $\_$ stands for the implicit index when the child of $v_s$ from previous step was promoted in $T$.
    \item If $d_l \geq d_s$, we proceed in the same way until we reach a vertex where this doesn't hold.
    \item Once we reach such vertex say $v_{q'}$, if $d_{q'}\leq d_q$, we return all the visited vertices and the corresponding indices as above. 
    \item If on the other hand $d_s>d_{q'}>d_q$, we construct $T''$ by swapping $v_q$ and $v_{q'}$ in $T'$. It is easy to observe that such tree will have smaller Sombor index than $T'$. We initialize the same algorithm on $T''$. 
\end{itemize}

It is clear that the algorithm above will terminate, giving us a sequence of vertices $u_1,\dots, u_m$ and indices $i_1,\dots, i_{m-1}$ with $d(u_m)<d(u_1)\leq d(u_2)\leq \dots \leq d(u_{m-1})$ and $i_1>i_2>\dots>i_{m-1}$. This pair of sequences indicates that the expression
$$E = \big(h_{i_1}(u_1) - h_{i_1}(u_2)\big) + \big(h_{i_2}(u_2) - h_{i_2}(u_3)\big) + \dots + \big(h_{i_{m-1}}(u_{m-1}) - h_{i_{m-1}}(u_m)\big)$$
is contained in ${\rm SO}(T_{r+i}) - {\rm SO}(T'_{r+i})$ (assuming $T'$ was adjusted as in the last step of the algorithm above whenever necessary).

We observe that this expression can be simplified using telescoping, since $h_a$ is increasing in $a$. Indeed

$$E \leq h_{i_1}(u_1) - h_{i_{m-1}}(u_m) < 0$$

Combining this with our induction hypothesis, it follows easily that 
$${\rm SO}(T_{r+i}) - {\rm SO}(T'_{r+i}) \leq 0$$
And by selection of $r$, we obtain the desired result directly.
\end{proof}

\acknowledgment{I would like to thank professor Slobodan Filipovski for introducing me to this problem and giving me valuable suggestions.}
\bibliographystyle{plain} 
\bibliography{biblio}

\begin{thebibliography}{10}

\bibitem{reti2021sombor}
A.~Ali, T.~Do{\v{s}}lic, and T.~R{\'e}ti.
\newblock On the sombor index of graphs.
\newblock {\em Contrib. Math}, 3:11--18, 2021.

\bibitem{ghanbari2021sombor}
S.~Alikhani and N.~Ghanbari.
\newblock Sombor index of certain graphs.
\newblock {\em arXiv preprint arXiv:2102.10409}, 2021.

\bibitem{das2021sombor}
I.~N. Cangul, A.~S. {\c{C}}evik, K.~C. Das, and Y.~Shang.
\newblock On sombor index.
\newblock {\em Symmetry}, 13(1):140, 2021.

\bibitem{liu2021more}
H.~Chen, X.~Fang, H.~Liu, Z.~Tang, and Q.~Xiao.
\newblock More on sombor indices of chemical graphs and their applications to
  the boiling point of benzenoid hydrocarbons.
\newblock {\em International Journal of Quantum Chemistry}, 121(17):e26689,
  2021.

\bibitem{cruz2021sombor}
R.~Cruz, I.~Gutman, and J.~Rada.
\newblock Sombor index of chemical graphs.
\newblock {\em Applied Mathematics and Computation}, 399:126018, 2021.

\bibitem{trees}
C.~Delorme, O.~Favaron, and D.~Rautenbach.
\newblock Closed formulas for the numbers of small independent sets and
  matchings and an extremal problem for trees.
\newblock {\em Discret. Appl. Math.}, 130(3):503--512, 2003.

\bibitem{gutman2021geometric}
I.~Gutman.
\newblock Geometric approach to degree-based topological indices: Sombor
  indices.
\newblock {\em MATCH Commun. Math. Comput. Chem}, 86(1):11--16, 2021.

\bibitem{gutman2021some}
I.~Gutman.
\newblock Some basic properties of sombor indices.
\newblock {\em Open J. Discret. Appl. Math}, 4(1):1--3, 2021.

\bibitem{liu2022sombor}
I~Gutman, Y.~Huang, H.~Liu, and L.~You.
\newblock Sombor index: Review of extremal results and bounds.
\newblock {\em Journal of Mathematical Chemistry}, pages 1--28, 2022.

\bibitem{milovanovic2021some}
M.~Matejic, E.~Milovanovic, and I.~Milovanovic.
\newblock On some mathematical properties of sombor indices.
\newblock {\em Bull. Int. Math. Virtual Inst}, 11(2):341--353, 2021.

\bibitem{rada2021general}
J.~Rada, J.~M Rodr{\'\i}guez, and J.~M Sigarreta.
\newblock General properties on sombor indices.
\newblock {\em Discrete Applied Mathematics}, 299:87--97, 2021.

\bibitem{redvzepovic2021chemical}
I.~Red{\v{z}}epovi{\'c}.
\newblock Chemical applicability of sombor indices.
\newblock {\em Journal of the Serbian Chemical Society}, 2021.

\bibitem{wang}
H.~Wang.
\newblock The extremal values of the wiener index of a tree with given degree
  sequence.
\newblock {\em Discrete Applied Mathematics}, 156, 10 2007.

\end{thebibliography}
\end{document}